\documentclass[12pt]{amsart}
\usepackage{amsfonts}

\usepackage[all,cmtip]{xy}
\usepackage{amsmath}
\usepackage{amsthm}
\usepackage{amssymb}
\usepackage{enumitem}

\newtheorem*{definition*}         {Definition}

\theoremstyle{remark}

\newcommand*{\Q}{\mathbb{Q}}

\newcommand*{\Z}{\mathbb{Z}}

\newcommand*{\GG}{\mathbb{G}}

\newcommand*{\C}{\mathbb{C}}

\usepackage{amscd,amssymb,amsmath}
\usepackage{color}
\definecolor{purple}{rgb}{0.59, 0.44, 0.84}

\begin{document}


%
%
\newcount\refCount
\def\newref#1 {\advance\refCount by 1
\expandafter\edef\csname#1\endcsname{\the\refCount}}
\newref AX
\newref AXB
\newref BERTRAND
\newref BERTRANDPILLAY
\newref BRION
\newref KIRBY
\newref OSWAL
\newref PT
\newref SIMPSONASPECTS
\newref SIMPSONSUBSPACES

%
\newcount\derandomCount
\def\derandom#1 {\advance\derandomCount by 1
\expandafter\edef\csname#1\endcsname{\the\derandomCount}}

\title{Ax-Schanuel and the rank one Riemann-Hilbert correspondence}

\author{Jonathan Pila and 
Jacob Tsimerman}

\maketitle

\section{Introduction}

Let $X$ be a smooth projective complex algebraic variety.
This note revisits a result of Simpson [\SIMPSONSUBSPACES, Theorem 3.1(c)] 
(reproved by Oswal [\OSWAL, Theorem 1.1] via o-minimal geometry) 
concerning subspaces that are (irreducible)
algebraic varieties on both sides of the rank one Riemann-Hilbert
correspondence for $X$: {\it such varieties must be cosets of algebraic subgroups.}
The result suggests an underlying Ax-Schanuel theorem.
We show that this is the case and 
recover the bi-algebraicity statement via functional transcendence.
We do not address Simpson's full result that bi-algebraic varieties 
are ``triple tori'' (i.e. they are also tori in the Dolbeault space), 
his characterization of the tori that are bi-algebraic, 
nor his description of those  defined over $\overline{\Q}$ 
on both sides, when $X/\overline{\Q}$.

The {\it rank $n$ Riemann-Hilbert correspondence\/} is a complex analytic isomorphism 
(biholomorphism)
$$
{\rm rihi}: M_{\rm dR}^n(X) \rightarrow M_B^n(X)
$$
between the {\it de Rham moduli space\/} $M_{\rm dR}^n(X)$ of rank $n$ vector bundles 
with flat connection on $X$ and the {\it Betti moduli space\/} $M_B^n(X)$ 
of rank $n$ complex local systems on $X$. 

In rank one, both source and target are 
commutative complex algebraic groups,
and ${\rm rihi}$ is a group isomorphism. On the de Rham side,
the connected component of the identity
$$
M_{\rm dR}^1(X)^{o}={\rm Pic}^0(X)^\natural,
$$
where $A^\natural$ denotes the universal vectorial extension of an abelian variety $A$ 
(see e.g. [\BERTRANDPILLAY, \BRION]),
while on the Betti side (with $g=\dim {\rm Pic}^0(X)$) we have
$$
M_B^1(X)^{o} = (\C^\times)^{2g}.
$$

Restriction to identity components yields a biholomorphic group 
isomorphism of (complex connected) commutative 
algebraic groups
$$
{\rm rihi}: {\rm Pic}^0(X)^\natural \rightarrow (\C^\times)^{2g}.
$$


Our Ax-Schanuel result concerns 
a biholomorphic group isomorphism
$$
\theta: R\rightarrow S
$$
of commutative algebraic groups (which are always assumed to be complex and connected
and are identified with their sets of complex points)
with a certain immiscibility property that we
formulate next. In fact, we do not need $\theta$ to be biholomorphic, only that its graph $\Theta$
is an analytic subgroup dominant to each factor with discrete fibres.
And we also need only a weaker immiscibility property: see \S6.



\medbreak
\noindent
{\bf 1.1. Definition.\/} 
A pair $R, S$ of commutative algebraic groups is called a
{\it disjoint pair\/} if every connected algebraic subgroup $W\subset R \times S$ is a product
$W=R'\times S'$, where $R'\subset R, S'\subset S$ are algebraic subgroups.

\medbreak
\noindent
{\bf 1.2. Definition.}
A commutative algebraic group $G$ is called {\it NVQ\/} (``no vectorial quotients'') 
if $G$ admits no non-trivial homomorphisms to a vector group.

\medbreak

Such groups are called {\it almost semi-abelian\/} in [\BERTRANDPILLAY], where
they also suggest the term {\it anti-additive} (in analogy with the notion
of anti-affine group [\BRION]).
Note that, for an abelian variety $A$, and for $h\ge 1$, 
the groups $A^\natural$ and $(\C^\times)^{h}$ are a disjoint pair, 
and both groups are NVQ (see [\BERTRANDPILLAY, \BRION]).

We defer the statement of our main result, Theorem 4.1, to \S4.
Here we state an implication of it in a special case that includes the 
rank one Riemann-Hilbert correspondence, and implies the bi-algebraicity statement;
it is an ``Ax-Lindemann'' theorem in this setting.


By a {\it coset\/} we will mean a coset of an irreducible algebraic subgroup.
We will say that a coset is NVQ if the corresponding subgroup is NVQ.
Subvarieties are assumed closed in their ambients.
The Zariski closure of $A\subset V$ is denoted ${\rm Zcl}(A)$; the relevant ambient
variety should always be clear from the context.

\medbreak
\noindent
{\bf 1.3. Theorem.\/} 
{\it Let $R, S$ be a disjoint pair of commutative algebraic groups and 
suppose that $\theta: R\rightarrow S$ is a biholomorphic isomorphism.
Let $V\subset R$ be an algebraic subvariety, $S'$ the smallest coset
of $S$ containing $\theta(V)$.
If $S'$ is NVQ then ${\rm Zcl}\big(\theta(V)\big)=S'$
(and so is a coset). \/}

\medbreak
\noindent
{\bf 1.4. Definition.\/}
By a {\it bi-algebraic subvariety\/} for $\theta: R\rightarrow S$ we
mean an (irreducible) algebraic subvariety $V\subset R$
whose $\theta$-image is an algebraic subvariety $W\subset S$.

\medbreak
\noindent
{\bf 1.5. Corollary.\/} {\it Let $R, S$ be a disjoint pair of commutative algebraic groups, 
and suppose that
$\theta: R\rightarrow S$ is a biholomorphic isomorphism. 
If $\theta:V\rightarrow W$ is bi-algebraic and the smallest coset of $S$
containing $W$ is NVQ, then $W$ and $V$ are cosets.\/} 

\medbreak

Corollary 1.5 recovers Simpson's bi-algebraicity result as the target group
is semi-abelian, and so all its cosets are NVQ.

\medbreak
\noindent
{\bf Proof of Corollary 1.5.\/} 
By Theorem 1.3, $W$ is a coset.
Then $V$ inherits a group structure
and so is also a coset.\ \qed

\medbreak

Note that the converse implication does not hold: 
cosets are in general not bi-algebraic (unless $\theta$ is algebraic, but this is
excluded by the disjoint pair assumption).

To prove our Ax-Schanuel theorem for $\theta: R\rightarrow S$ we rely 
on an Ax-Schanuel theorem for commutative algebraic groups established in [\PT], 
which generalizes an Ax-Schanuel theorem
for NVQ commutative algebraic groups due to Ax-Bertrand [\AXB, \BERTRAND]. 

Finally we note that our statement (though adequate for deducing Theorem 1.3)
seems to be weaker than the optimal statement one would like;
see \S4.

\medskip
\noindent
{\bf Ackowledgements.\/} The authors are 
very grateful the Institute for Advanced 
Study for its support and warm hospitality 
in the period that this paper was written.

\medskip
\noindent
{\bf AI disclosure.\/} No AI was used 
in the research or writing of this 
paper.

\medbreak

\section{Ax-Schanuel for a commutative algebraic group}

\medskip

A commutative algebraic group $G$ has a canonical connected
subgroup $G_{\rm NVQ}$ which is NVQ and such that 
$$
G \cong G_{\rm NVQ} \times G_{\rm VEC}
$$
for some (non-canonical) vector group $G_{\rm VEC}$ (see [\PT, Lemma 2.1]). 
The group $G_{\rm NVQ}$ will be called the {\it maximal NVQ subgroup\/} of $G$.

For a commutative algebraic group $G$ we write
$$
{\rm nvqdim}(G)=\dim (G_{\rm NVQ}).
$$

Note that, for commutative algebraic groups $G_1,\ldots, G_k$ with $H_i$ the maximal
NVQ subgroup of $G_i$, 
$$
\big(G_1\times\ldots\times G_k\big)_{\rm NVQ}= H_1\times\ldots\times H_k,
$$
and so 
$$
{\rm nvqdim}\big(G_1\times\ldots\times G_k\big)=
{\rm nvqdim} H_1 + \ldots + {\rm nvqdim} H_k.
$$

We can then state Ax-Schanuel for a commutative algebraic group as in [\PT].
Note that [\PT, Cor. 2.3] is stated erroneously in the original arXiv version, 
corrected in Corollary 2.2 below.

The universal cover of $G$ is $\tilde{G}=\C^d$, where $d=\dim G$, 
and the universal covering map $\pi: \tilde{G}\rightarrow G$ 
is invariant under some discrete group of translations.

\medbreak
\noindent
{\bf 2.1. Theorem.\/} ([\PT, Thm. 2.2]) 
{\it Let $H$ be an NVQ commutative algebraic group, and $V$ a vector group. 
Let $\pi: \tilde{H}\rightarrow H$ be the universal covering map of $H$ 
with graph $\Gamma=\Gamma_\pi$.
Let $Z\subset (\tilde{H}\times H)\times V$ be an algebraic variety
and 
$$
U\subset Z\cap (\Gamma \times V)
$$ 
a complex-analytically irreducible component. Then
$$
\dim Z \ge \dim U +\dim H
$$
unless the image of $U$ under projection to $H$ lies in a strict coset.\ \qed
}

\medbreak

The theorem may be applied to a commutative algebraic group $G$ putting $H=G_{\rm NVQ}$
and with $V=G_{\rm VEC}$ any choice of vectorial complement.
Though the decomposition of $G$ is non-canonical,
one can show that, if the image of $U$ under projection lies in a strict coset,
then the connected subgroup underlying the strict coset is independent
of the choice of vectorial complement.

\medbreak
\noindent
{\bf 2.2. Corollary.} ([\PT, Cor. 2.3])
{\it
Let $H$ be an NVQ commutative algebraic group, and $V$ a vector group. 
Let $\pi: \tilde{H}\rightarrow H$ be the universal cover
with graph $\Gamma=\Gamma_\pi$.
Let $Z\subset (\tilde{H}\times H) \times V$ be an algebraic variety
and 
$$
U\subset Z\cap(\Gamma \times V)
$$ 
a complex-analytically irreducible component. 
Let $J$ be a minimal coset containing the image of $U$ 
under projection to $H$. Then
$$
\dim Z \ge \dim U + {\rm nvqdim} (J). \ \qed
$$
}





\bigbreak

\section{Ax-Schanuel in a differential field}

\medskip

We work in a differential field 
$(K, \{D_1,\ldots, D_k\}, \C)$ where 
$\Q\subset K$, the $D_i$ are
commuting derivations, and the constant field is $\C$.

Let $G$ be a commutative algebraic group. We take some embedding of $G$
as a quasi-projective variety. Under this embedding, the
$K$-points of $G$ are denoted $G(K)$.

If $\underline{x}=(x_1,\ldots, x_n)\in G(K)$
then a {\it lift\/} $\tilde{\underline{x}}\in \tilde{G}(K)$ 
of $\underline{x}$ is a tuple of $K$-points (in $K^{\dim G}$)
characterized by a suitable algebraic differential equation. 
In the classical case of the exponential function 
$G=(\C^\times)^n$ these
are the differential equations $D_i\tilde{x_j}=D_ix_j/x_j$ for each derivation $D_i$
and coordinate $x_j$; Ax-Schanuel here is proved by Ax [\AX].
In general, the differential equation is that of the logarithmic derivative 
mapping $G$ to its Lie algebra
and can be expressed in terms of invariant differentials on $G$.
If 
$$
\sum R_{ij}(\underline{X})dX_n, \quad i=1,\ldots, \dim G
$$ 
are $\C$-linearly independent 
invariant differentials then the equations take the form
$$
D_k\tilde{X}_j =\sum R_{ij}(\underline{X}) D_k X_j, \quad j=1,\ldots, d, \quad {\rm all\ } i {\rm\ and\ } k.
$$
The solutions are determined up to translation by constants.
(See e.g. Kirby [\KIRBY] in the semi-abelian case; Bertrand [\BERTRAND] in general.)

The {\it rank\/} of $\underline{x}\in K^n$, denoted ${\rm rank}(\underline{x})$, 
is the $K$-rank of the matrix $(D_ix_j)$.
Set ${\rm trd}(\underline{x})= {\rm trd}_\C\C(\underline{x})$.

If $G$ is a commutative algebraic group then there is an algebraic isomorphism (over $\C$)
$$
G\cong G_{\rm NVQ}\times G_{\rm VEC}
$$
for some (non-canonical) 
vector group $G_{\rm VEC}$.
Any such isomorphism induces an isomorphism
$$
G(K) \rightarrow G_{\rm NVQ}(K)\times G_{\rm VEC}(K),
$$
under which the image of $\underline{x}$ will be denoted
$$
\underline{x}\mapsto (\underline{x}_{\rm N}, \underline{x}_{\rm V}).
$$

Now we can state Ax-Schanuel for a commutative algebraic group $G$ in the 
setting of a differential field.

\medbreak
\noindent
{\bf 3.1. Theorem.\/} 
{\it Let $G$ be a commutative algebraic group 
and $\underline{x}\in G(K)$. Fix a decomposition 
$G\cong G_{\rm NVQ}\times G_{\rm VEC}$.
Let $\underline{\tilde{x}}_{\rm N}\in \tilde{G}_{\rm NVQ}(K)$ 
be a lift of $\underline{x}_{\rm N}$.
Then
$$
{\rm trd} (\underline{x}, \underline{\tilde{x}}_{\rm N}) \ge 
{\rm rank} (\underline{x})+ {\rm nvqdim}(G)
$$
unless $\underline{x}_{\rm N}\in J(K)$ for some strict coset $J\subset G_{\rm NVQ}$.
}

\medbreak
\noindent
{\bf Proof.\/} This follows from Theorem 2.1 via Seidenberg embedding. \qed

\medbreak
\noindent
{\bf 3.2. Corollary.\/}
{\it 
Let $G$ be a commutative algebraic group and $\underline{x}\in G(K)$. 
Fix a decomposition 
$G\cong G_{\rm NVQ}\times G_{\rm VEC}$.
Let $\underline{\tilde{x}}_{\rm N}\in \tilde{G}_{\rm NVQ}(K)$ 
be a lift of $\underline{x}_{\rm N}$.
Then
$$
{\rm trd} (\underline{x}, \underline{\tilde{x}}_{\rm N}) \ge 
{\rm rank} (\underline{x})+ {\rm nvqdim}(J)
$$
where $J\subset G_{\rm NVQ}$ is a minimal coset such that $\underline{x}_{\rm N}\in J(K)$.
}

\medbreak
\noindent
{\bf Proof.\/} This follows likewise from Corollary 2.2. \qed

\bigbreak

\section{Main Theorem}

\medskip

In this section we consider a biholomorphic isomorphism $\theta: R\rightarrow S$
of a disjoint pair of commutative algebraic groups. (Hence the minimal coset appearing 
in the statement of Theorem 4.1 is a product.)

Let $\Theta\subset R\times S$ be the graph of $\theta$. The (any) lift $\tilde{\Theta}$ is a linear
subvariety, hence $\Theta$ is a differential algebraic subgroup and we can consider
its $K$-points.

For a coset $R'\times S'$ we define 
$$
\delta(R'\times S')=\max\big({\rm nvqdim} R', {\rm nvqdim} S'\big).
$$

Our main theorem is the following.

\medbreak
\noindent
{\bf 4.1. Theorem.\/} (Weak Ax-Schanuel for a biholomorphic disjoint pair)
{\it 
Let $\theta: R\rightarrow S$ is a biholomorphic isomorphism of a disjoint pair
of commutative algebraic groups.
Suppose $(\underline{x}, \underline{y})\in \Theta(K)$. Let $R'\times S'$ be a minimal coset
containing the image of $(\underline{x}, \underline{y})$ under projection to
$R_{\rm NVQ}\times S_{\rm NVQ}$. Then
$$
{\rm trd} (\underline{x}, \underline{y}) \ge {\rm rank\/}(\underline{x}, \underline{y})
+ \delta(R'\times S').
$$

}

\medbreak
\noindent
{\bf Proof.\/} We have a decomposition $G\cong G_{\rm NVQ}\times G_{\rm VEC}$ inducing
a complex algebraic isomorphism of their $K$-points, as above. So we have
$$
\C(\underline{x})= \C(\underline{x}_N, \underline{x}_V), \quad
\C(\underline{y})= \C(\underline{y}_N, \underline{y}_V).
$$
By the assumptions, 
$$
\underline{x}_N\in R'(K),\quad \underline{y}_N \in S'(K).
$$
We take some decompositions
$$
R' \cong R'_{\rm NVQ} \times R'_{\rm VEC},\quad 
S' \cong S'_{\rm NVQ} \times S'_{\rm VEC}.
$$
These give isomorphisms of their $K$-points, under which
$$
\underline{x}_N \mapsto (\underline{x}_{NN}, \underline{x}_{NV}), \quad 
\underline{y}_N \mapsto (\underline{y}_{NN}, \underline{y}_{NV}),
$$
so that
$$
\C(\underline{x}_N)=\C(\underline{x}_{NN}, \underline{x}_{NV}), \quad
\C(\underline{y}_N)=\C(\underline{y}_{NN}, \underline{y}_{NV}).
$$
By Corollary 3.2 we have
$$
{\rm trd}(\underline{x}, \tilde{\underline{x}}_N, 
\underline{y}, \tilde{\underline{y}}_N) \ge 
{\rm rank} (\underline{x}, \underline{y}) + 
{\rm nvqdim}(R')+ {\rm nvqdim}(S').
$$

To complete the proof we require the upper estimate 
$$
{\rm trd} \big(\C(\underline{x}, \tilde{\underline{x}}_N, \underline{y}, 
\tilde{\underline{y}}_N)/
\C(\underline{x}, \underline{y})\big)\le
\min\big({\rm nvqdim} R', {\rm nvqdim} S'\big).
$$
We first show that the transcendence degree is at most ${\rm nvqdim}(R')$.

The isomorphism $\theta$ lifts to a linear isomorphism 
of the covering spaces. 
Thus there is a $\C$-affine relationship between any pair of $K$-points 
$(\tilde{\underline{x}}, \tilde{\underline{y}})\in \tilde{\Theta}(K)$ lifting 
$(\underline{x}, \underline{y})\in \Theta(K)$, i.e. 
$$
\tilde{\theta}(\tilde{\underline{x}}) = \tilde{\underline{y}} + \underline{c}
$$ 
for some constant complex vector $\underline{c}$.
Thus
$$
\C(\underline{x}, \tilde{\underline{x}}, \underline{y}, 
\tilde{\underline{y}})=
\C(\underline{x}, \tilde{\underline{x}}, \underline{y})
$$
and
$$
{\rm trd} \big(\C(\underline{x}, \tilde{\underline{x}}, \underline{y}, 
\tilde{\underline{y}})/
\C(\underline{x}, \underline{y})\big)\le
{\rm trd}\big(\C(\underline{x}, \tilde{\underline{x}})/
\C(\underline{x})\big).
$$
Moreover, there are $\C$-linear isomorphisms on the covering spaces 
$$
\tilde{R} \cong \tilde{R}_{N} \times \tilde{R}_V, \quad
\tilde{R'} \cong \tilde{R'}_N \times \tilde{R'}_V
$$
giving $\C$-affine isomorphisms on $K$-points
(as we are assuming $(\underline{x}_N\in R'(K)$)
under which
$$
\tilde{\underline{x}} \mapsto (\tilde{\underline{x}}_N, 
\tilde{\underline{x}}_{V}) 
\mapsto 
(\tilde{\underline{x}}_{NN}, \tilde{\underline{x}}_{NV}, \tilde{x}_V).
$$

As $R_V$ and $R'_V$ are vector groups they are isomorphic 
as algebraic groups with
their universal covers giving, isomorphisms on $K$-points
under which, up to translations,
$$
\tilde{\underline{x}}_{V}\mapsto \underline{x}_V, \quad 
\tilde{\underline{x}}_{NV}\mapsto \underline{x}_{NV}.
$$

Thus
$$
\C(\underline{x}, \tilde{\underline{x}}_N, \underline{y}, 
\tilde{\underline{y}}_N)=
\C(\underline{x}, \tilde{\underline{x}}, \underline{y})=
\C(\underline{x}, \tilde{\underline{x}}_{NV}, \underline{y})
$$
and
$$
{\rm trd}\big(\C(\underline{x}, \tilde{\underline{x}}_N, \underline{y}, 
\tilde{\underline{y}}_N)/
\C(\underline{x}, \underline{y})\big) \le
{\rm trd}\big(\C(\underline{x}, 
\tilde{\underline{x}}_{NV})/\C(\underline{x})\big)
\le {\rm nvqdim}(R').
$$
Likewise,
$$
{\rm trd} \big(\C(\underline{x}, \tilde{\underline{x}}_N, \underline{y}, 
\tilde{\underline{y}}_N)/
\C(\underline{x}, \underline{y})\big)
\le {\rm nvqdim}(S').
$$
by a symmetrical argument on the $S$-side. \qed

\medskip

An adequate Ax-Schanuel theorem should be sharp in the sense that the 
transcendence degree estimate should be an equality for the full locus 
(i.e. a $K$-point of full rank)
for the differential-algebraic varieties governing the 
asserted transcendence degree,
in this case the loci of the form $\Theta\cap(R'\times S')$.
The above might not be always sharp, so we call it ``weak''. 

However, it is sharp for the ambient $R\times S$ for a disjoint pair $R, S$ of 
NVQ commutative algebraic groups. 
For then the differential-algebraic
variety $\Theta$ has 
$$
{\rm rank}(\Theta) = d=\dim R =\dim S,
$$
while $\tilde{\Theta}$ is linear, 
and for a $K$-point 
$(\underline{x}, \underline{y})\in\Theta(K)$
of rank $d$ we have
$$
{\rm trd}(\underline{x}, \underline{y})= {\rm rank}(\underline{x}, 
\underline{y}) + 
\dim(R\times S) - \dim \tilde{\Theta} = 2d.
$$

\medbreak

\noindent
{\bf 4.2. Corollary.\/} {\it With the hypotheses of the theorem, 
assume further
that $R, S$ are NVQ commutative algebraic groups of dimension $d$,
and that $(\underline{x}, \underline{y})$ is not contained in any strict
coset.
Then
$$
{\rm trd}(\underline{x}, \underline{y}) \ge 
{\rm rank}(\underline{x}, \underline{y})+d.
\qed
$$
\/}

\medbreak

We state a complex-analytic version of Theorem 4.1.

\medbreak
\noindent{\bf 4.3. Theorem.\/}
{\it 
Let $\theta: R\rightarrow S$ is a biholomorphic isomorphism of a disjoint pair
of commutative algebraic groups.
Let $\Theta\subset R\times S$ be the graph of $\theta$,
$Z\subset R\times S$ an algebraic subvariety, $U$
a complex-analytically irreducible component of $Z\cap\Theta$,
and $R'\times S'$ the minimal coset containing the image of $U$ under
projection to $R_{\rm NVQ}\times S_{\rm NVQ}$. Then
$$
\dim Z \ge \dim U + \delta(R'\times S'). \qed
$$

}

\medskip

Here we would expect that a ``sharp'' Ax-Schanuel would assert some equality
(with the ``right'' $\delta$), as the complex
analytic component $U$ cannot have smaller dimension than the conditions imposed
by the intersection inside a suitable ambient. But here we don't
know what the right ambient or $\delta$ is.
Thus here again the conclusion may not be sharp, but is so
in the setting corresponding to Corollary 4.2.

\medbreak
\noindent
{\bf 4.4. Corollary.\/}
{\it
With the hypotheses of Theorem 4.3 and assuming that $R, S$ are
NVQ commutative algebraic groups of dimension $d$,
$$
\dim Z =\dim U + d
$$
unless $U$ is contained in a strict coset.\/}

\bigbreak

\section{Proof of Theorem 1.3}

\medskip
\noindent
{\bf Proof of Theorem 1.3.\/} 
Let (as exists in some suitable differential field $K$)
$$
(\underline{x}, \underline{y})\in \big(\Theta \cap (V \times S)\big)(K)
$$
be of full rank 
$$
d={\rm rank}(\underline{x}, \underline{y})=\dim V.
$$

Let $W={\rm Zcl}(\theta(V))$, so that $\underline{y}\in W(K)$.
Thus $S'$ is the smallest coset containing $W$.
Let $R''$ be the smallest coset of $R$ containing $V$, so that
$$
R'' \times S'
$$
is the smallest coset in $R\times S$ such that
$(\underline{x}, \underline{y})\in (R''\times S')(K)$.

We apply Theorem 4.1 but note that the target group $S'$
in the statement is NVQ, so $S'_{\rm NVQ}=S$. 
We fix some decomposition
$$
R'' \cong R_{\rm NVQ} \times R_{\rm VEC}
$$
and observe that the restriction 
$$
R'' \times S' \rightarrow R_{\rm NVQ} \times S_{\rm NVQ}
$$
of the projection $R\times S \rightarrow R_{\rm NVQ} \times S_{\rm NVQ}$
is the identity on the second factor.

Let $R'$ be the smallest coset containing the image of $V$ under
projection to $R_{\rm NVQ}$, so that
$$
R'\times S'
$$
is the minimal coset containing the image of $(\underline{x}, \underline{y})$
under projection to $R_{\rm NVQ} \times S_{\rm NVQ}$
(this is still a coset in $R\times S$ and so is a product by the
disjoint pair assumption).

By Theorem 4.1 we have
$$
{\rm trd} (\underline{x}, \underline{y}) \ge 
{\rm rank} (\underline{x}, \underline{y})
+ \dim S'.
$$
Therefore we have
$$
d+\dim S' \le {\rm trd} (\underline{x}, \underline{y}) \le d +\dim W
$$
and since $W\subset S'$ we conclude $W=S'$.\ \qed

\medbreak

We can analyse bi-algebraic subvarieties a bit further in the case that
$\theta: A^\natural\rightarrow (\C^\times)^{2g}$, where $A$
is an abelian variety with $\dim A=g$. Suppose $\theta(V)=W$ is bi-algebraic,
so that $W\cong (\C^\times)^f$. On the other side, $V$ is some additive
extension of an abelian subvariety $B\subset A$, and considering
the rank of the deck group we must have $\dim B=h$ and $\dim V=f=2h$.
All such extensions occurring as cosets in $A^\natural$ are non-trivial, 
so $V=B^\natural$.

\medbreak

\section{A variant}

\medskip

Here we want to consider biholomorphic isormorphisms 
$\theta : A^\natural \rightarrow B^\natural$ where
$A, B$ are abelian varieties of dimension $g$. Following the notation
and characterization in [\BRION], we denote by
$$
A^\natural_{\rm Aff} \cong \C^g
$$
the {\it affine subgroup\/}. 
It is the kernel of the extension $A^\natural \rightarrow A$,
and, when $A$ is simple, it is the unique strict
positive dimensional connected algebraic subgroup of $A^\natural$
(see [\BRION, Thm. 3.3]).

Note that no such pair is a disjoint pair, as there are (an abundance of) non-product 
subgroups in $A^\natural_{\rm Aff} \times B^\natural_{\rm Aff}$.

However, the disjoint property is more than we need to prove
4.1: we need it only
for the subgroups that arise as intersections with $\Theta$.

\medbreak
\noindent{\bf 6.1. Definition.\/} Let $\theta: R\rightarrow S$ be
a biholomorphic group isomorphism of commutative algebraic groups (with graph $\Theta$).
We say that $R, S$ are $\Theta$-disjoint if, for every  
algebraic subgroup $W\subset R\times S$, and connected component $\Phi\subset W\cap\Theta$,
the smallest coset containing $\Phi$ is a product $R'\times S'$ of cosets.

\medbreak

The proof of Theorem 4.1 only requires product structure on the cosets arising after 
intersection with $\Theta$, because the $K$-point 
$(\underline{x}, \underline{y})\in \Theta(K)$. 
Thus, the conclusions of 4.1, 4.2, (4.3, 4.4,) 1.3, and 1.5
hold under the weaker asumption of $\Theta$-disjointness.

\medbreak
\noindent
{\bf 6.2. Proposition.\/}
{\it Suppose $\theta(A^\natural_{\rm Aff})=B^\natural_{\rm Aff}$.
Then $\theta$ descends to an algebraic isomorphism $A\rightarrow B$.
\/}

\medbreak
\noindent
{\bf Proof.\/}
We have the projections $\pi_A: A^\natural\rightarrow A$ and 
$\pi_B: B^\natural\rightarrow B$. 
Suppose $\pi_A(x)=\pi_A(x')$. Then $\pi_A(x-x')=0_A$ so
$x-x'\in A^\natural_{\rm Aff}$ and $\theta(x-x')\in B^\natural_{\rm Aff}$.
And then $\pi_B(\theta(x))=\pi_B(\theta(x'))$.
So $\theta$ induces a complex analytic map $A\rightarrow B$,
algebraic by Chow's theorem.\ \qed

\medbreak

At the other extreme, 
restricting to the case that $A, B$ are simple and
non-isogenous, if
$\theta(A^\natural_{\rm Aff}) \cap B^\natural_{\rm Aff}=\{(0,0)\}$
we get $\Theta$-disjoint pairs.

\medskip
\noindent
{\bf 6.3. Proposition.\/}
{\it Let $A, B$ be simple and non-isogenous abelian varieties
of dimension $g$. If $W\subset A^\natural\times B^\natural$ is a non-product coset
then $W\subset A^\natural_{\rm Aff} \times B^\natural_{\rm Aff}$.
\/}

\medbreak
\noindent
{\bf Proof.\/} We will make heavy use of a theorem of Brion 
[\BRION, Prop. 3.3]: every strict subgroup scheme $H\subset A^\natural$
is affine (and likewise for $B^\natural$).

Suppose $W\subset A^\natural \times B^\natural$ is a
connected algebraic subgroup.
Let
$$
W_x=\{y\in B^\natural: (x,y) \in W\}, \quad x\in A^\natural.
$$
Then $W_x$ is a coset of the subgroup $W_0$, 
which is either $B^\natural$ (leading to a product)
or (we may assume) a finite union of cosets of
subgroups of $B^\natural_{\rm Aff}$ over a finite subgroup
of $B$.

Suppose $W$ is not dominant to $B^\natural$. Then it sits over
a finite subgroup of $B$, and so sits over $B^\natural_{\rm Aff}$
by its connectedness. So $W\subset A^\natural \times B^\natural_{\rm Aff}$.
Then we consider fibres $W^y$ over $y\in B^\natural$.
If $W$ is not dominant to $A^\natural$ then again it lies over
$A^\natural_{\rm Aff}$ (by connectedness) and so 
$W\subset A^\natural_{\rm Aff}\times B^\natural_{\rm Aff}$.

So assume $W$ is dominant to $A^\natural$. If the fibres $W_y=A^\natural$ 
then $W$ is a product.
If not, the fibres are affine groups over some finite subgroup of $A$. 
If $W$ is dominant to $A^\natural$,
the map sending $y\in \pi_B(W)$ to the sum of the projections of 
the fibres is a dominant map to $A$,
which is impossible.

So we can suppose $W$ is dominant to $B^\natural$ and then the 
sum of the projections of the fibres
is a dominant map from $A^\natural\rightarrow B$. But this map is 
constant on the (affine) fibres
of $A^\natural \rightarrow A$, and so we have $A\rightarrow B$
dominant,
which contradicts our hypothesis.

Thus, the only non-product cosets are from inside 
$A^\natural_{\rm Aff} \times B^\natural_{\rm Aff}$. \qed

\medbreak
\noindent{\bf 6.4. Corollary.\/}
{\it If $A, B$ are simple and non-isogenous, and
$\theta(A^\natural_{\rm Aff}) \cap B^\natural_{\rm Aff}=\{(0,0)\}$,
then $A^\natural, B^\natural$ are $\Theta$-disjoint.\ \qed}

\medbreak
\noindent{\bf 6.5. Proposition.\/}
{\it In the setting of Corollary 6.4 the only bi-algebraic variety
of positive dimension is $A^\natural$.}

\medbreak
\noindent
{\bf Proof.\/}
Let $V\subset A^\natural$ of positive dimension.
If $V\subset A^\natural_{\rm Aff}$ then
$\theta(V)\not\subset B^\natural_{\rm Aff}$ and 
so the smallest coset containing $\theta(V)$ is $B^\natural$,
which is NVQ. Then ${\rm Zcl}(\theta(V))=B^\natural$
by the conclusion of 1.5 under $\Theta$-disjointness.
A symmetric argument applies if $W\subset B^\natural$ of positive dimension.
If $\Theta: V\rightarrow W$ is bi-algebraic of positive dimension, then either $V$ or $W$ is not
contained in the respective affine part, and the Zariski closure of the image
of (respectively) $\theta^{-1}(W)$ or $\theta(V)$ is (respectively) $A^\natural$ or $B^\natural$.
And so $V=A^\natural$ or $W=B^\natural$, and so both hold. \qed

\bigbreak
\bigbreak

\centerline{References}

\medskip

\begin{enumerate}

\item[\AX.] J. Ax, On Schanuel’s conjectures, 
{\it Annals Math. (2)\/} {\bf 93} (1971), 252–268.

\item[\AXB.] J. Ax, 
 Some topics in differential algebraic geometry I : Analytic subgroups of
algebraic groups,
{\it Amer. J. Math.} {\bf  94} (1972), 1195-1204.

\item[\BERTRAND.]
D. Bertrand, 
Schanuel’s conjecture
for non-isoconstant elliptic curves over function fields,
{\it Model theory with applications to algebra and analysis. Vol. 1.
(Newton Institute, Cambridge, 2005)\/}, 41-62.
Z. Chatzidakis, D. Macpherson, A. Pillay, and A. Wilkie, editors,
LMS Lecture Notes {\bf 349}, CUP, Cambridge, 2008.

\item[\BERTRANDPILLAY.]
D. Bertrand and A. Pillay,
A Lindemann-Weierstrass theorem for semi-abelian varieties 
over function fields,
{\it J. Amer. Math. Soc.\/} {\bf 23} (2010), 491–533.

\item[\BRION.] M. Brion,
Anti-affine algebraic groups,
{\it Journal of Algebra\/} {\bf 321\/} (2009), 934–952.



\item[\KIRBY.] J. Kirby,
The theory of the exponential differential equations of semiabelian varieties,
{\it Selecta Math. (N.S.)} {\bf 15\/} (2009), 445–486.

\item[\OSWAL.] A. Oswal, 
Bi-algebraicity in the rank one Riemann-Hilbert correspondence 
via o-minimality,
{\it Advances in Mathematics\/} {\bf 448} (2024), 109704.

\item[\PT.] J. Pila and J. Tsimerman,
Ax-Schanuel and exceptional integrability, arXiv:2202.04023.

\item[\SIMPSONASPECTS.] C. Simpson,
Transcendental aspects of the Riemann-Hilbert correspondence,
{\it Illinois J. Math\/} {\bf 34} (1990),  368--391.

\item[\SIMPSONSUBSPACES.] C. Simpson,
Subspaces of moduli spaces of rank one local systems,
{\it Annales Scientifiques de l’\'E.N.S. (4)} {\bf 26} (1993), 361-401.



\end{enumerate}

\bigbreak
\bigbreak
\bigbreak

\rightline{20261006}

\bigbreak
\bigbreak
\bigbreak

\medskip

\end{document}